\documentclass[11pt]{article}
\usepackage[utf8]{inputenc}
\usepackage[english]{babel}
\usepackage{amsmath}
\usepackage{amsfonts}
\usepackage{amssymb}
\usepackage{graphicx}
\usepackage{pgf,tikz}
\usepackage{mathrsfs}
\usetikzlibrary{arrows}

\def\R{\mathbb R}

\def\N{\mathbb N}

\author{G.~Chadzitaskos, 
M.~Havl\' \i  \v cek, J.~Patera\\
}
\title{Orthonormal bases on $L^2(\R^+)$}

\begin{document}

\maketitle
\abstract{We derive the explicit form of eigenvectors of selfadjoint extension $H_\xi$, parametrized by 
$\xi \in \langle 0,\pi),$ of differential expression $ H=-\frac{d^2 }{d x^2} + \frac{x^2 }{4}$ together with the spectrum $\sigma(H_\xi)$  on the space $L^2(\R^+).$  For each $\xi$ the set of eigenvectors  form an orthonormal basis of $L^2(\R^+).$

\section{Introduction}
The basic examples of quantum mechanics is a quantization  the harmonic oscillator. A selfadjoint Hamiltonian  $H_D$ of the one--dimensional linear harmonic oscillator is  generated by the differential expression
 \begin{equation}{\label{DifOp}}
 H=-\frac{d^2  }{d x^2} + \frac{x^2 }{4} 
  \end{equation}
  with appropriate definition domain $D.$ It is known that the operator $H_D$ has a pure point spectrum and its eigenfunctions form the orthonormal basis in $L^2(\R),$ and $H_D$ is a unique selfadjoint operator generated by $H$ on $L^2(\R).$
  
  The situation is quite different on $L^2(\R^+),$ there is one--parametric set of  selfadjoint operators $H_\xi,$ $\xi\in\langle 0, \pi)$ with corresponding definition domains  $D_\xi$ and with the same differential expression (\ref{DifOp}) (\cite{BEH},p. 137). All these selfadjoint operators are selfadjoint extensions of the closed symmetric operator $\hat{H}$ with the domain $\hat{D}= \bigcap\limits_{\xi\in\langle 0, \pi)} D_\xi.$ Following the theorem (\cite{W} p. 246) all these extension have the same essential spectrum. As in the case of the operator $H_{\xi=0}, $ where it applies $\sigma_{ess} (H_{\xi=0}) = \emptyset,$  it applies for all operators $H_\xi,$ $\xi\in\langle 0, \pi).$ In other words, for any $\xi\in\langle 0, \pi)$ there exist an orthonormal basis formed by eigenvectors of $H_\xi.$ The objective of this paper is to derive  explicit form of the orthonormal basis and express $\sigma(H_\xi).$

\section{Parabolic cylinder functions}

 The parabolic cylinder fuctions \cite{RG} (9.240, 9.210)
\begin{equation}{\label{PCF}}  
   D_{\nu}(x) = e^{-\frac{x^2}{4}} [\frac{\sqrt{\pi} \, 2^{\frac{\nu}{2}} }{\Gamma (\frac{1-\nu}{2})}  {_1} \Phi_1(-\frac{\nu}{2},\frac{1}{2},\frac{x^2}{2})- \frac{\sqrt{\pi} \, 2^{\frac{\nu+1}{2}} }{\Gamma (\frac{-\nu}{2})}\, x \,{_1} \Phi_1(\frac{1-\nu}{2},\frac{3}{2},\frac{x^2}{2})]
    \end{equation} 
   
  are the solutions of the Weber differential equation \cite{RG} (9.255) 
  $$(\frac{d^2  }{d x^2} - \frac{x^2 }{4} + \nu +\frac{1}{2})D_\nu (x) =0.$$
% that approach zero for $x\rightarrow \infty.$ 
Values $\nu \in \lbrace 0, 1, 2, \dots \rbrace\equiv \N_0 $ need special attention, because of
$$\frac{1}{\Gamma (\frac{-\nu}{2})}=0,\,\,\, \Gamma (\frac{1-\nu}{2})= \infty\,\, \mbox{ for } \nu = 1,3,5,\ldots $$
and
$$\frac{1}{\Gamma (\frac{-\nu}{2})}=\infty,\,\,\, \Gamma (\frac{1-\nu}{2})= 0\,\, \mbox{ for } \nu = 0,2,4,\ldots .$$

Definition (\ref{PCF}) then gives $D_\nu (x) = h_\nu (x)$ known hermitian functions \cite{RG} (9.253).

The following relations  holds for PCFs \cite{RG}(7.711, 8.370, 8.372):
\begin{equation}{\label{Norm}}
\int_0^\infty |D_\nu(x)|^2 dx = \frac{1}{c(\nu)^2},\,\,\,c(\nu)=\sqrt{\sqrt{\frac{2}{\pi}}\frac{\Gamma(-\nu)}{ \beta(-\nu)}},\,\,\, \beta(-\nu)=\sum_{k=0}^\infty \frac{(-1)^k}{-\nu + k}
\end{equation}
(note that $c(\nu) \, \, D_\nu$ is normalized), and

\begin{equation}{\label{Ss}}
\int_0^\infty D_\nu(x) D_\mu(x) dx = \frac{\pi \, 2^{\frac{1}{2} (\nu +\mu +1)}}{\mu - \nu} [\frac{1}{\Gamma(\frac{1-\mu}{2})\Gamma(-\frac{\nu}{2})} - \frac{1}{\Gamma(\frac{1-\nu}{2})\Gamma(-\frac{\mu}{2})}].
\end{equation}

\section{Two Lemmas and two Theorems}
It is known that the differential expression (\ref{DifOp})
$$
 H=-\frac{d^2  }{d x^2} + \frac{x^2 }{4} 
$$
  with definition domain 
  \begin{equation}
 D_\xi\,(H): =  \lbrace f \in\tilde{ \cal D} , f(0) \cos \xi - f'(0) \sin \xi =0 \rbrace ,
\end{equation}
  is a selfadjoint operator on $ L^2(R^+)$ for all $\xi \in \langle 0,\pi),$  and $\tilde{ \cal D}= \lbrace f \, f' \in a.c.(0,\infty): f, H\,f, \in L^2(R^+) \rbrace $ \cite{BEH} (p. 127, p. 137)

So, if $D_\nu$  will belong to ${\cal D_\xi}(H)$ for some $\nu$  then $D_\nu$ will be an eigenvector of the considered selfadjoint operator with eigenvalues $\nu+1/2$. Eq. (\ref{Norm}) guarantees that $D_\nu$ lies in $L^2 (\R^+)$).   

The last condition generates the relationship
\begin{equation}{\label{00}}
 D_\nu(0) \cos \xi - D_\nu'(0) \sin \xi =0. 
\end{equation}

Although, values $D_\nu(0)$ and $D_\nu'(0)$ can be calculated using  definition (\ref{PCF}), we have to distinguish two cases:
% and eq. 4 has the following form ($\xi \neq 0$  )
\begin{enumerate}
\item when $\nu \notin \N_0$ we  obtain
\begin{equation}{\label{Eqv}}
\eta \Gamma(-\frac{\nu}{2})-\Gamma(\frac{1-\nu}{2})=0, \, \, \eta=\frac{1}{\sqrt{2}} \cot \xi.
\end{equation}
\item  when $\nu  \in \N_0$ we obtain
\begin{equation}{\label{Herm00}}
 h_\nu(0) \cos \xi - h_\nu'(0) \sin \xi =0. 
\end{equation} 
If  $\nu$ is odd, then $ h_\nu(0)=0,\,\,\, h_\nu'(0) =1,$ and Eq. (\ref{Herm00}) is fulfilled only if $\xi=0.$ If $\nu$ is even, then  $ h_\nu(0)=1,\,\,\, h_\nu'(0) =0,$ and Eq. (\ref{Herm00}) is fulfilled only if $\xi=\frac{\pi}{2}.$  In both cases, condition (\ref{Herm00}) is fulfilled by the set of Hermitian functions $\lbrace h_0, h_2,h_4,\ldots \rbrace $ and $\lbrace h_1, h_3,h_5,\ldots \rbrace ,$ respectively. It is known that both sets form orthonormal bases in $L^2(\R^+).$ 

\end{enumerate}

Eq.(\ref{Eqv}) has to be solved for $\nu.$

First we prove two lemmas.
\vspace{.2 cm}

{\bf Lemma 1}:
\begin{enumerate}

\item If $\nu  \in (2M-1,2M), M=1,2,\ldots$ or $\nu < 0,$ then $\beta(-\nu) \geq 0,$ 

\item If $\nu  \in (2M-2,2M-1), M=1,2,\ldots,$  then $\beta(-\nu) < 0.$ 

\end{enumerate}

{\bf Proof:}

1. Using the relationship
$$\beta(-\nu)=\sum_{k=0}^\infty \frac{1}{(-\nu + 2k)(-\nu + 2k+1)}, $$ 
\cite{RG}(8.372), it is possible to show by elementary calculation that $(-\nu + 2k)(-\nu + 2k+1)>0 $ for all $k=0,1,\ldots $ if $\nu \in (2M-1,2M), M=0,1,\ldots,$ or $\nu<0.$

2. In this case we rewrite the sum $\beta(-\nu)$ in the following form:

$$\beta(-\nu)= -\frac{1}{\nu}+\sum_{k=0}^\infty \frac{1}{-\nu + k+1} +\frac{1}{-\nu + k+2}= \\
-\frac{1}{\nu}-\sum_{k=0}^\infty\frac{1}{(-\nu + k+1)(-\nu + k+2)}. $$

For considered values of $\nu \in (2M-2,2M-1), M=1,2,\ldots $ the products $(-\nu + k+1)(-\nu + k+2)$ are positive, and therefore all denominators of the members in the previous sum are positive and so $\beta \leq 0$ (note that $-\frac{1}{\nu}<0 $). $\square$

Remark: Comparing functions $\Gamma(-\nu)$ and $\beta(-\nu)$ we have the relationship
$$\mbox{sgn} ( \Gamma(-\nu))= \mbox{sgn} (\beta(-\nu)),\, \nu \in \R. $$

It shows that normalization factor $c(\nu)$ (Eq.\ref{Norm})) is correctly defined. 
\vspace{.2 cm}
{\bf Lemma 2}: 

Function $y(\nu):= \frac{\Gamma(\frac{1-\nu}{2})}{\Gamma(-\frac{\nu}{2})}$ has the following properties:

\begin{enumerate} 

\item
There are asymptotes for $\nu_{\mbox{as}} \in \lbrace 2n+1| n \in \N_0 \rbrace$ and $\lim_{\nu\rightarrow\nu_{\mbox{as}}^+} y(\nu)=\infty,$ $\lim_{\nu\rightarrow\nu_{\mbox{as}}^-} y(\nu)=-\infty.$ Further $\lim_{\nu\rightarrow - \infty} y(\nu)=+\infty.$
\item The set  $ \lbrace 2n| n \in \N_0 \rbrace$ consists of all zero points of $y$.
\item In the intervals  $(-\infty,1),(M,M+1), M=0,1,\ldots, $ $y$ is continuous decreasing function.
\end{enumerate}

{\bf Proof:}

1. The first assertion is a direct consequence of explicit form \cite{RG}(8.314) of function $\Gamma.$ 

For the remaining assertions it is sufficient to prove that the sequence $ \lbrace y(-2n)| n \in \N_0 \rbrace$ is growing and $\lim_{n\rightarrow  \infty} y( - 2 n)=+\infty.$ As  $y( - 2 n)$ can be expressed as

$$  y( - 2 n)= \frac{\Gamma(n+\frac{1}{2})}{\Gamma(n)} =  \frac{\sqrt{ \pi}(2n-1)!!}{2^{n}(n-1)!} $$ \cite{RG}(8.339), the assertion can be easily verified.

2. $\Gamma$ -- function has no zero points. Therefore $y(\nu)=0$ only if $|\Gamma(-\frac{\nu}{2})|=\infty,$ i. e. $\nu=2n.$

3. For $y'(\nu)$ we obtain
$$\frac{dy(\nu)}{d\nu}=\frac{1}{2} \frac{\Gamma(\frac{1-\nu}{2})}{\Gamma(-\frac{\nu}{2})} \, [\psi(-\frac{\nu}{2})-\psi(\frac{1-\nu}{2})],\,\,\psi(\mu) = \frac{d}{d \mu}  \lg \Gamma(\mu). $$

Using the relationships
$$\psi(-\frac{\nu}{2})-\psi(\frac{1-\nu}{2})=-2\, \beta(-\nu), \mbox{ and } \Gamma(-\nu) =  \frac{2^{-\nu-1}}{\sqrt{\pi}}\Gamma(\frac{1-\nu}{2})\,\Gamma(-\frac{\nu}{2}), $$
\cite{RG}(8.370,8.335) we obtain
$$\frac{dy(\nu)}{d\nu} =  -2^{\nu}\frac{\sqrt{\pi}} {\Gamma(-\frac{\nu}{2})^2}\Gamma(-\nu) \,\beta (-\nu).$$
As $ \Gamma(-\nu) \,\beta (-\nu)>0$ (see Remark) the proof is completed. $\square$

The consequence of this Lemma is a Theorem

{\bf Theorem 1:}

  For any $\eta \in \R$ and any $M \in \N=\lbrace1,2,\ldots\rbrace$ there is just one solution $\nu_\eta^{(M)}$ of  Eq. (\ref{Eqv})  in the interval $I_M,$ where
 $$I_1= (-\infty,1), \,\, I_M= (2M-1, 2M+1),\, M=2,3,\ldots$$

No further solution of Eq. (\ref{Eqv}) exists.

\begin{table}[h]
\caption{Example of first 11 values of $\nu_\eta$}
\begin{tabular}{|l|l|l|l|l|l|l|}

\hline
 $\nu_{-2.18}$& $\nu_{-0.51}$& $\nu_{0}$ &$\nu_{ 0.23}$& $\nu_{0.51}$&$\nu_{0.97}$&$\nu_{2.18}$  \\
\hline
0.77051& 0.399912& 0&-0.311391& -0.875066&-2.33401& -9.95 \\
\hline
 2.66471& 2.26065& 2.&1.86885& 1.71369& 1.5141&1.26337 \\
\hline
 4.59639& 4.20523& 4.&3.90249& 3.78578&  3.62177&3.36297 \\
\hline 
 6.54652& 6.1743& 6.& 5.91892& 5.82117& 5.67849& 5.42659\\
\hline 
 8.50776& 8.15402& 8.& 7.92911& 7.84326& 7.715227& 7.47292\\
\hline 
 10.4764& 10.1394& 10.& 9.93622& 9.85874& 9.74156&9.50897\\
\hline
 12.4503& 12.1283& 12.& 11.9415& 11.8704&11.7617& 11.5382\\
 \hline
 14.4281& 14.1195& 14.&  13.9457&13.8795& 13.7777&13.5626\\
\hline 
 16.409& 16.1123& 16.&  15.9491&15.887& 15.7908&15.5834 \\
\hline
 18.3922& 18.1062& 18.& 17.9519& 17.8932&17.8019& 17.6014  \\
 \hline
20.3773& 20.101& 20.&  19.9543& 19.8985& 19.8113&19.6172  \\
\hline
\end{tabular}
\end{table}
\begin{figure}[h]
\includegraphics{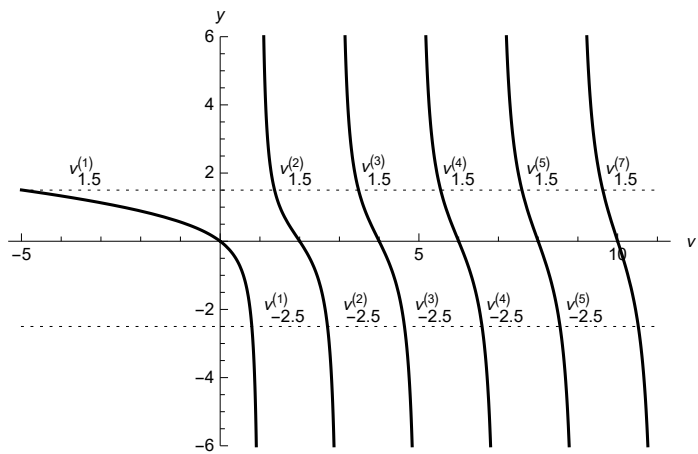}
\caption{Function $y(\nu) = \frac{ \Gamma \left(\frac{1-\nu}{2}\right)}{\Gamma \left(-\frac{\nu}{2}\right)}$}
\end{figure}

Let $\Omega_\xi$ denote  the set

$$\Omega_\xi=\lbrace \nu^{(M)}_{\cot \xi},
\,  M=1,2,\ldots\rbrace ,\, \xi \in \langle 0,\pi),
\,  \nu^{(M)}_{\cot \xi}= \nu^{(M)}_{\eta}, $$  (we understand $\Omega_0 = \lbrace 0, 2, 4, \ldots\rbrace),$ 

and let denote further by $\cal{E}_\xi$ the set

$${\cal E}_\xi =
\lbrace c (\nu) D_{\nu} |\,\nu \in \Omega_\xi\rbrace.$$

The set $\cal{E}_\xi \subset D_\xi$ contains all eigenvectors of the selfadjoint operator $H_\xi,$ and the set $ \lbrace \nu+\frac{1}{2}, \nu \in \Omega_\xi\rbrace$ contains all eigenvalues of  $H_\xi.$ 

Note that orthogonality of two different eigenvectors can bee seen also from the Eq.(\ref{Ss}). For different $\mu, \,\,\nu$ fulfilling the condition  
$\Gamma(\frac{1-\mu}{2})/ \Gamma(-\frac{\mu}{2}))=\Gamma(\frac{1-\nu}{2})/\Gamma(-\frac{\nu}{2})=\eta,$ which is our case, Eq. (\ref{Ss}) is the scalar product in $L^2(R^+) $ equal to zero. Moreover, the Eq. (\ref{Norm}) guarantees that the eigenvectors are normalized.

We denote further by $\hat{H}$ the restriction of $H_\xi$
to domain
$$\hat{D} =\lbrace f \in \tilde{D},\, f(0)=f'(0)=0 \rbrace \subset D_\xi (H) \subset L^2(\R^+). $$

Operator $\hat{H}$ is closed, symmetric with deficiency indices (1,1) \cite{BEH} (prop. 4.8.6, p.129), and
$H_\xi$ is a selfadjoint extension of $\hat{H}$
for any $\xi \in \langle 0, \pi).$ Selfadjoint extensions $H_{\xi=0}$ and  $H_{\xi=\frac{\pi}{2}}$ have pure point spectra, which is equivalent to the existence of orthonormal bases in $L^2(\R^+).$
The basis is in the case  $D_{\xi=0}(H)=\lbrace h_{2n+1}|n=0,1,2,\ldots\rbrace, $ and it is in the case  $D_{\xi=\frac{\pi}{2}}(H)=\lbrace h_{2n+1}|n=0,1,2,\ldots\rbrace .$ 
As we mentioned in the introduction, the same  is true for all operators $H_\xi$ with any parameter $\xi.$

Consequently, one can write theorem

{ \bf Theorem 2}

The set $\cal{E}_\xi  $ consisting of eigenvectors of $H_\xi$ is an orthonormal  basis in $L^2(\R^+)$ for any $\xi \in \langle 0, \pi),$ 
and $\sigma(H_\xi)= \lbrace \nu+\frac{1}{2}, \nu \in \Omega_\xi\rbrace.$

\section{Concluding remarks}

The results we present can be translated to the case $L^2(\R^-).$ In this case  an orthonormal basis in $L ^2(\R^-)$ is
$$\tilde{\cal{E}_\xi}: = \lbrace
\tilde{D}_{\nu}|\, \nu \in \Omega_\xi  \rbrace,\tilde{D}_\nu(x):= D_\nu(-x).
$$

 These two bases can be combined to the base in $L^2(\R).$ As $L^2(\R)=L^2(\R^+) \oplus L^2(\R^-).$ Then for any pair $(\xi, \sigma) \in \langle 0, \pi) \times  \langle 0, \pi)  $    the set $\cal{E}_\xi \oplus \tilde{\cal{E}_\sigma}$ is a basis in $L^2(\R).$ Explicitly
$${\cal{E}_\xi} \oplus \tilde{\cal{E}}_\sigma = \lbrace
(D_\nu, 0),\, \nu \in \Omega _\xi  \rbrace \cup\lbrace
(0, \tilde{D}_\nu),\,  \nu \in \Omega _\sigma  \rbrace .
$$
Of note, the known orthonormal basis $\lbrace  h_n, n=0,1,\ldots \rbrace$ of $ L^2(\R)$ consisting of hermitian functions is not contained in this set. Functions $h_n$ are eigenvectors of selfadjoint operator $H_D$ with definition domain
$$D=\lbrace f, f' \mbox{absolutely continuous}, f ,\,Hf \in L^2(\R)\rbrace, $$
and operator $H_D$ is physically interpreted as Hamiltonian of quantum linear harmonic oscillator.

\section*{Acknowledgement}
All authors acknowledge the financial support from RVO14000 and "Centre for Advanced Applied Sciences", Registry No. CZ.02.1.01/0.0/0.0/16\_019/0000778, supported by the Operational Programme Research, Development and Education, co-financed by the European Structural and Investment Funds and the state budget of the Czech Republic.

The authors (G.C., M.H.) thank P. Exner for helpfull discussion.

 \end{document}